\documentclass{amsproc}%
\usepackage{graphicx}
\usepackage{amscd}
\usepackage{amsmath}
\usepackage{amsfonts}
\usepackage{amssymb}%
\theoremstyle{plain}

\numberwithin{equation}{section}

\begin{document}
\title[INTEGRALITY OVER\ FIXED\ RINGS OF\ AUTOMORPHISMS]{INTEGRALITY OVER\ FIXED\ RINGS OF\ AUTOMORPHISMS IN\ A\ LIE\ NILPOTENT\ SETTING}
\author{Jen\H{o} Szigeti}
\address{Institute of Mathematics, University of Miskolc, Miskolc, Hungary 3515}
\email{matjeno@uni-miskolc.hu}
\thanks{The author was partially supported by the National Research, Development and Innovation}
\thanks{Office of Hungary (NKFIH) K119934.}

\begin{abstract}
Let $R$ be a Lie nilpotent algebra of index $k\geq1$ over a field $K$ of
characteristic zero. If $G$ is an $n$-element subgroup $G\subseteq
\mathrm{Aut}_{K}(R)$ of $K$-automorphisms, then we prove that $R$ is right
integral over $\mathrm{Fix}(G)$ of degree $n^{k}$. In the presence of a
primitive $n$-th root of unity $e\in K$, for a $K$-automorphism $\delta
\in\mathrm{Aut}_{K}(R)$ with $\delta^{n}=\mathrm{id}_{R}$, we prove that the
skew polynomial algebra $R[w,\delta]$ is right integral of degree $n^{k}$ over
$\mathrm{Fix}(\delta)[w^{n}]$.

\end{abstract}
\maketitle

\noindent1. INTRODUCTION

\bigskip

The notions of algebraicity and integrality play an important role in algebra
and algebraic number theory. The aim of our paper is to present two
integrality results, which are non-commutative generalizations of well known
facts, the first one is closely related to classical Galois theory and the
second one is related to the integrality of the commutative polynomial algebra
$C[w]$ over $C[w^{n}]$.

A natural replacement of commutativity is the following: a ring $R$\ is called
Lie nilpotent of index $k$, if the left normed commutator product (with
$[x,y]=xy-yx$)%
\[
\lbrack\lbrack\ldots\lbrack\lbrack x_{1},x_{2}],x_{3}],\ldots,x_{k}%
],x_{k+1}]=0
\]
is a polynomial identity on $R$.

In section (5) we prove the following two integrality theorems.

\bigskip

\noindent\textbf{Theorem (A).}\textit{ Let }$R$\textit{ be a Lie nilpotent
algebra of index }$k\geq1$\textit{ over a field }$K$\textit{ of characteristic
zero. If }$G$\textit{ is an }$n$\textit{-element subgroup }$G\subseteq
\mathrm{Aut}_{K}(R)$\textit{ of }$K$\textit{-automorphisms, then }$R$\textit{
is right integral over }$\mathrm{Fix}(G)$\textit{ of degree }$n^{k}$\textit{.
In other words, for any }$r\in R$\textit{ we have}%
\[
c_{0}+rc_{1}+\cdots+r^{n^{k}-1}c_{n^{k}-1}+r^{n^{k}}=0
\]
\textit{for some }$c_{i}\in\mathrm{Fix}(G)$\textit{, }$0\leq i\leq n^{k}%
-1$\textit{.}

\bigskip

\noindent\textbf{Remarks to (A).} Let $K\subseteq L$ be a field extension and
$G\subseteq\mathrm{Aut}_{K}(L)$ be a finite subgroup of the $K$-automorphisms
of $L$. Then $a\in L$ is a root of the monic polynomial $f(w)=\Pi_{\gamma\in
G}(w-\gamma(a))$ of degree $n=\left\vert G\right\vert $ and the coefficients
of $f(w)\in L[w]$\ are in the fixed intermediate field $K\subseteq
\mathrm{Fix}(G)\subseteq L$\ of $G$. Thus $L$ is integral over $\mathrm{Fix}%
(G)$ of degree $n$. We note that the replacement of $L$ by a commutative
$K$-algebra $R$ (in the above argument) immediately yields Theorem (A) for
$k=1$.

\noindent If $\mathrm{Fix}(G)=K$ and $a\in L$ is also a root of a
$K$-irreducible polynomial $q(w)\in K[w]$, then $q(w)$ is a divisor of $f(w)$.
For a finite normal field extension $K\subseteq L$, this observation shows
that the Galois group $G=\mathrm{Aut}_{K}(L)$ acts transitively on the roots
of $q(w)$.

\bigskip

\noindent\textbf{Theorem (B).}\textit{ Let }$R$\textit{ be a Lie nilpotent
algebra of index }$k\geq1$\textit{ over a field }$K$\textit{ of characteristic
zero. If }$e\in K$\textit{ is a primitive }$n$\textit{-th root of unity and
}$\delta\in\mathrm{Aut}_{K}(R)$\textit{ is a }$K$\textit{-automorphism with
}$\delta^{n}=\mathrm{id}_{R}$\textit{, then the skew polynomial algebra
}$R[w,\delta]$\textit{ is right integral of degree }$n^{k}$\textit{ over the
}$K$\textit{-subalgebra }$\mathrm{Fix}(\delta)[w^{n}]$\textit{ of }%
$R[w,\delta]$\textit{. In other words, for any }$f(w)\in R[w,\delta]$\textit{
we have}%
\[
g_{0}(w^{n})+f(w)g_{1}(w^{n})+\cdots+f^{n^{k}-1}(w)g_{n^{k}-1}(w^{n}%
)+f^{n^{k}}(w)=0
\]
\textit{for some }$g_{i}(w^{n})\in\mathrm{Fix}(\delta)[w^{n}]$\textit{,
}$0\leq i\leq n^{k}-1$\textit{ (notice, that }$w^{n}$\textit{ is central in
}$R[w,\delta]$\textit{).}

\bigskip

\noindent\textbf{Remarks to (B).} For a commutative $K$-algebra $C$, the
choice $R=C$ and $\delta=\mathrm{id}_{C}$ in Theorem (B) immediately yields
the well-known fact, that the polynomial algebra $C[w]$ is integral of degree
$n$ over the $C$-subalgebra $C[w^{n}]$ generated by the power $w^{n}$. A
simple illustration of the case $n=2$ is the following: any $f(w)\in C[w]$ can
be uniquely written as $f(w)=p(w^{2})+wq(w^{2})$ and $f^{2}(w)-2p(w^{2}%
)f(w)+(p^{2}(w^{2})-w^{2}q^{2}(w^{2}))=0$ proves the integrality of $C[w]$
over $C[w^{2}]$ of degree $2$.

\bigskip

The proofs of Theorem A and B are heavily based on the use of certain matrix
algebras over the given Lie nilpotent algebra $R$. The so called Lie nilpotent
Cayley-Hamilton theorem from [S1] is an indispensable ingredient of our
development. In order to provide a self-contained treatment, we present the
necessary prerequisites in sections (2),(3) and (4).

Throughout the paper an algebra $R$\ means a not necessarily commutative
unitary algebra over a field $K$ of characteristic zero, all $K$-subalgebras
contain the identity and all $K$-endomorphisms preserve the identity. The
group of units in $R$ is denoted by $\mathrm{U}(R)$ and the centre of $R$ is
denoted by $\mathrm{Z}(R)$. The notation for the full $n\times n$ matrix
algebra over $R$ is $\mathrm{M}_{n}(R)$ and $\mathrm{GL}_{n}(R)=\mathrm{U}%
(\mathrm{M}_{n}(R))$. The matrix $E_{i,j}\in\mathrm{M}_{n}(R)$ has $1$ in the
$(i,j)$ position and zeros in all other positions. The fixed $K$-subalgebra of
a $K$-endomorphism $\delta:R\longrightarrow R$ is $\mathrm{Fix}(\delta)=\{r\in
R\mid\delta(r)=r\}$.

\bigskip

\noindent2. THE\ ALGEBRA OF\ SKEW CENTRALIZING\ MATRICES

\bigskip

Using the natural (element-wise) extension $\delta_{n}:\mathrm{M}%
_{n}(R)\longrightarrow\mathrm{M}_{n}(R)$ of a $K$-endomorphism $\delta
:R\longrightarrow R$ and a matrix $W\in\mathrm{M}_{n}(R)$, one can define the
subalgebra of skew $(\delta,W)$-centralizing matrices in $\mathrm{M}_{n}%
(R)$\ as%
\[
\mathrm{M}_{n}(R,\delta,W)=\{A\in\mathrm{M}_{n}(R)\mid\delta_{n}(A)W=WA\}.
\]
For the sake of brevity, in the rest of the paper we omit the qualifier "skew".

If $W\in\mathrm{M}_{n}(\mathrm{Fix}(\delta))$, then $\mathrm{M}_{n}%
(R,\delta,W)$ is closed with respect to the action of$~\delta_{n}$.

\bigskip

\noindent\textbf{Example (1).} The natural $\mathbb{Z}_{2}$-grading
$E=E_{0}\oplus E_{1}$ defines an automorphism $\varepsilon(g_{0}+g_{1}%
)=g_{0}-g_{1}$ of the Grassmann algebra
\[
E=K\left\langle v_{1},\ldots,v_{i},\ldots\mid v_{i}v_{j}+v_{j}v_{i}=0\text{
for all }1\leq i\leq j\right\rangle
\]
generated by the infinite sequence of anticommutative indeterminates
$(v_{i})_{i\geq1}$

\noindent with $v_{i}^{2}=0$.

\noindent Using the invertible diagonal matrix%
\[
Q_{d}=E_{1,1}+\cdots+E_{d,d}-E_{d+1,d+1}-\cdots-E_{n,n}\text{ }%
\]
in $\mathrm{M}_{n}(K)$\textit{, }we obtain the classical supermatrix algebra
$\mathrm{M}_{n,d}(E)$ as%
\[
\mathrm{M}_{n,d}(E)=\mathrm{M}_{n}(E,\varepsilon,Q_{d})=\{A\in\mathrm{M}%
_{n}(E)\mid\varepsilon_{n}(A)Q_{d}=Q_{d}A\}.
\]
The shape of a matrix $A\in\mathrm{M}_{n,d}(E)$ is%
\[
A=\left[
\begin{array}
[c]{cc}%
A_{1,1} & A_{1,2}\\
A_{2,1} & A_{2,2}%
\end{array}
\right]  ,
\]
where the square blocks $A_{1,1}$ and $A_{2,2}$ are of sizes $d\times d$ and
$(n-d)\times(n-d)$ and the rectangular blocks $A_{1,2}$ and $A_{2,1}$ are of
sizes $d\times(n-d)$ and $(n-d)\times d$. The entries of $A_{1,1}$ and
$A_{2,2}$ are in the even part $E_{0}$ of $E$, while the entries of $A_{1,2}$
and $A_{2,1}$ are in the odd part $E_{1}$ of $E$.

\noindent We note that $E$ is Lie nilpotent of index $2$, the $K$-algebras
$\mathrm{M}_{n}(E)$ and $\mathrm{M}_{n,d}(E)$ play an important role in
Kemer's classification of T-prime T-ideals (see [Ke]). $\square$

\bigskip

\noindent\textbf{Example (2).} A $\mathbb{Z}_{n}$-grading of $R$ is an
$n$-tuple $(R_{0},R_{1},\ldots,R_{n-1})$, where each $R_{i}$ is a $K$-subspace
of $R$ such that
\[
R=R_{0}\oplus R_{1}\oplus\cdots\oplus R_{n-1}%
\]
is a direct sum and $R_{i}R_{j}\subseteq R_{i+j}$ for all $i,j\in
\{0,1,\ldots,n-1\}$, where $i+j$ is taken in $\{0,1,\ldots,n-1\}$\ modulo $n$.
Using such a $\mathbb{Z}_{n}$-grading of $R$, the $K$-subalgebra
$\mathrm{M}_{n}^{g}(R)$ of $\mathrm{M}_{n}(R)$ is defined as follows:%
\[
\mathrm{M}_{n}^{g}(R)=\{A=[a_{i,j}]\in\mathrm{M}_{n}(R)\mid a_{i,j}\in
R_{j-i}\text{ for all }1\leq i,j\leq n\}=
\]%
\[
\left[
\begin{array}
[c]{ccccc}%
R_{0} & R_{1} & \cdots & R_{n-2} & R_{n-1}\\
R_{n-1} & R_{0} & R_{1} & \ddots & R_{n-2}\\
\vdots & R_{n-1} & \ddots & \ddots & \vdots\\
R_{2} & \ddots & \ddots & R_{0} & R_{1}\\
R_{1} & R_{2} & \cdots & R_{n-1} & R_{0}%
\end{array}
\right]  ,
\]
where $j-i$ is taken in $\{0,1,\ldots,n-1\}$\ modulo $n$.

\noindent If $e\in K$ is a primitive $n$-th root of unity, then%
\[
\widehat{e}(r_{0}+r_{1}+\cdots+r_{n-1})=r_{0}+er_{1}+\cdots+e^{n-1}r_{n-1}%
\]
defines a $K$-automorphism $\widehat{e}:R\longrightarrow R$ of the
$\mathbb{Z}_{n}$-graded $K$-algebra $R=R_{0}\oplus R_{1}\oplus\cdots\oplus
R_{n-1}$ and%
\[
D_{e}=e^{-1}E_{1,1}+e^{-2}E_{2,2}+\cdots+e^{-n}E_{n,n}%
\]
is a diagonal matrix in $\mathrm{GL}_{n}(K)$\textit{.}\ Now $\mathrm{M}%
_{n}^{g}(R)$\ can be obtained as the algebra of $(\widehat{e},D_{e}%
)$-centralizing matrices:%
\[
\mathrm{M}_{n}^{g}(R)=\mathrm{M}_{n}(R,\widehat{e},D_{e})=\{A\in\mathrm{M}%
_{n}(R)\mid\widehat{e}_{n}(A)D_{e}=D_{e}A\}.\text{ }\square
\]

\bigskip

One of the remarkable properties of $(\delta,W)$-centralizing matrices is the following.

\bigskip

\noindent\textbf{2.1. Proposition.}\textit{ If }$W\!\in\!\mathrm{GL}%
_{n}(R)\cap\mathrm{M}_{n}(\mathrm{Z}(R))$\textit{ is an invertible matrix with
central elements, then the trace of a matrix }$A\!\in\!\mathrm{M}_{n}%
(R,\delta,W)$\textit{ is in the fixed ring of }$\delta$\textit{, i.e.
}$\mathrm{tr}(A)\in\mathrm{Fix}(\delta)$\textit{.}

\bigskip

\noindent\textbf{Proof.} Since the entries of $W$ are central and
$W^{-1}W=I_{n}$, a straightforward computation gives that $\mathrm{tr}%
(WAW^{-1})=\mathrm{tr}(A)$ for all $A\in\mathrm{M}_{n}(R)$. Thus for a matrix
$A\in\mathrm{M}_{n}(R,\delta,W)$ we have%
\[
\delta(\mathrm{tr}(A))=\mathrm{tr}(\delta_{n}(A))=\mathrm{tr}(WAW^{-1}%
)=\mathrm{tr}(A),
\]
whence $\mathrm{tr}(A)\in\mathrm{Fix}(\delta)$ follows. $\square$

\bigskip

Any $K$-endomorphism $\delta:R\longrightarrow R$ can be naturally extended to
a $K$-endomor-

\noindent phism $\delta_{z}:R[z]\longrightarrow R[z]$\ of the polynomial ring
$R[z]$: for $r_{0},r_{1}\ldots,r_{t}\in R$ take%
\[
\delta_{z}(r_{0}+r_{1}z+\cdots+r_{t}z^{t})=\delta(r_{0})+\delta(r_{1}%
)z+\cdots+\delta(r_{t})z^{t}.
\]

The elements of the skew polynomial ring $R[w,\delta]$ in the skew
indeterminate $w$ are left polynomials of the form $f(w)=r_{0}+r_{1}%
w+\cdots+r_{t}w^{t}$ with $r_{0},r_{1},\ldots,r_{t}\in R$. Besides the obvious
addition, we have the following multiplication rule in $R[w,\delta]$:
$wr=\delta(r)w$. If $\delta^{n}=\mathrm{id}_{R}$ (such a $\delta$\ is an
automorphism), then $w^{n}$ is a central element of $R[w,\delta]$.

Let $G=\{\sigma_{1},\sigma_{2},\ldots,\sigma_{n}\}\subseteq\mathrm{Aut}%
_{K}(R)$ be an $n$-element subgroup of the group of all $K$-automorphisms of
$R$ (with $\sigma_{n}=\mathrm{id}_{R}$). The fixed $K$-subalgebra of $G$ is
\[
\mathrm{Fix}(G)=\{r\in R\mid\sigma(r)=r\text{ for all }\sigma\in G\}=%
{\textstyle\bigcap\nolimits_{\sigma\in G}}
\mathrm{Fix}(\sigma).
\]

Any element $\tau\in G$\ defines a permutation $\pi\in\mathrm{S}_{n}$ by%
\[
\tau\circ\sigma_{1}=\sigma_{\pi(1)},\ldots,\tau\circ\sigma_{n}=\sigma_{\pi(n)}%
\]
and the corresponding $n\times n$\ permutation matrix is%
\[
P_{\tau}=E_{1,\pi(1)}+\cdots+E_{n,\pi(n)}.
\]

\bigskip

\noindent\textbf{2.2. Theorem.}\textit{ Let }$G=\{\sigma_{1},\sigma_{2}%
,\ldots,\sigma_{n}\}\subseteq\mathrm{Aut}_{K}(R)$\textit{ be an }%
$n$\textit{-element subgroup of }$K$\textit{-automorphisms. If }$r\in
R$\textit{, then}%
\[
\Gamma_{G}(r)=\underset{i=1}{\overset{n}{%
{\textstyle\sum}
}}\sigma_{i}(r)E_{i,i}%
\]
\textit{defines a natural diagonal embedding }$\Gamma_{G}:R\longrightarrow
\mathrm{M}_{n}(R)$\textit{ of rings. For }$r\in R$\textit{ and }%
$c\in\mathrm{Fix}(G)$\textit{ we have }$\Gamma_{G}(cr)=c\Gamma_{G}%
(r)$\textit{, }$\Gamma_{G}(rc)=\Gamma_{G}(r)c$\textit{ and }$\Gamma
_{G}(c)=cI_{n}$\textit{.}

\noindent\textit{If }$\tau\in G$\textit{, then }$\Gamma_{G}:R\longrightarrow
\mathrm{M}_{n}(R,\tau,P_{\tau})$\textit{\ is an embedding of }$K$%
\textit{-algebras, where }$P_{\tau}\in\mathrm{GL}_{n}(K)\mathit{\ }$\textit{is
the permutation matrix defined by }$\tau$\textit{.}

\bigskip

\noindent\textbf{Proof.} Since $\Gamma_{G}(r)$ is diagonal, the verification
of the fact that $\Gamma_{G}$ is an injective ring homomorphism is straightforward.

\noindent If $r\in R$ and $c\in\mathrm{Fix}(G)$, then $\Gamma_{G}%
(cr)=c\Gamma_{G}(r)$, $\Gamma_{G}(rc)=\Gamma_{G}(r)c$ and $\Gamma
_{G}(c)=cI_{n}$ follow from the fact that $\sigma(cr)=c\sigma(r)$,
$\sigma(rc)=\sigma(r)c$ and $\sigma(c)=c$ for all $\sigma\in G$.

\noindent The calculations%
\[
\tau_{n}(\Gamma_{G}(r))=\underset{i=1}{\overset{n}{%
{\textstyle\sum}
}}\tau(\sigma_{i}(r))E_{i,i}=\underset{i=1}{\overset{n}{%
{\textstyle\sum}
}}\sigma_{\pi(i)}(r)E_{i,i}%
\]
and%
\[
\tau_{n}(\Gamma_{G}(r))P_{\tau}=\left(  \Sigma_{i=1}^{n}\sigma_{\pi
(i)}(r)E_{i,i}\right)  \left(  \Sigma_{j=1}^{n}E_{j,\pi(j)}\right)
=\Sigma_{i=1}^{n}\sigma_{\pi(i)}(r)E_{i,i}E_{i,\pi(i)}=
\]%
\[
\Sigma_{i=1}^{n}\sigma_{\pi(i)}(r)E_{i,\pi(i)}=\Sigma_{i=1}^{n}E_{i,\pi
(i)}\sigma_{\pi(i)}(r)E_{\pi(i),\pi(i)}=
\]%
\[
\left(  \Sigma_{i=1}^{n}E_{i,\pi(i)}\right)  \left(  \Sigma_{k=1}^{n}%
\sigma_{k}(r)E_{k,k}\right)  =P_{\tau}\Gamma_{G}(r)
\]
show that $\Gamma_{G}(r)\in\mathrm{M}_{n}(R,\tau,P_{\tau})$. $\square$

\bigskip

\noindent\textbf{2.3. Theorem.}\textit{ If }$\delta:R\longrightarrow
R$\textit{\ is a }$K$\textit{-endomorphism and }$W\in\mathrm{M}_{n}%
(R)$\textit{, then any homomorphism }$\varphi:R\longrightarrow\mathrm{M}%
_{n}(R,\delta,W)$\textit{\ of }$K$\textit{-algebras\ with }$\varphi\circ
\delta=\delta_{n}\circ\varphi$\textit{\ has a unique }$K$\textit{-homomorphic
extension}%
\[
\varphi^{(\delta)}:R[w,\delta]\longrightarrow\mathrm{M}_{n}(R[z])
\]
\textit{such that }$\varphi^{(\delta)}(r)=\varphi(r)$\textit{ and }%
$\varphi^{(\delta)}(w)=Wz$\textit{. For }$\varphi^{(\delta)}$\textit{ we have}%
\[
\varphi^{(\delta)}(r_{0}+r_{1}w+\cdots+r_{t}w^{t})=\varphi(r_{0}%
)+\varphi(r_{1})Wz+\cdots+\varphi(r_{t})W^{t}z^{t}.
\]

\noindent(1)\textit{ If }$f(w)\in R[w,\delta]$\textit{ and }$\varphi
(c)=cI_{n}$\textit{ for some }$c\in R$\textit{, then }$\varphi^{(\delta
)}(cf(w))=$

\noindent\ $\ \ \ \ c\varphi^{(\delta)}(f(w))$\textit{.}

\noindent(2)\textit{ If }$\varphi$\textit{ is injective and }$W\in
\mathrm{GL}_{n}(R)$\textit{ is invertible, then }$\varphi^{(\delta)}%
$\textit{\ is also injective.}

\noindent(3)\textit{ If }$W\in\mathrm{M}_{n}(\mathrm{Fix}(\delta))$\textit{,
then }$\varphi^{(\delta)}$\textit{ is an }$R[w,\delta]\longrightarrow
\mathrm{M}_{n}(R[z],\delta_{z},W)$\textit{ homomor-}

\noindent\ \ \ \ \ \textit{phism of }$K$\textit{-algebras.}

\bigskip

\noindent\textbf{Proof.} The $K$-linearity of $\varphi^{(\delta)}$\ is clear.
In order to see the multiplicative property of $\varphi^{(\delta)}$, it is
enough to prove that%
\[
\varphi^{(\delta)}((rw^{i})(sw^{j}))=\varphi^{(\delta)}(rw^{i})\varphi
^{(\delta)}(sw^{j})
\]
for all $r,s\in R$ and $0\leq i,j$.

\noindent Using $\varphi\circ\delta=\delta_{n}\circ\varphi$ and $\varphi
(s)\in\mathrm{M}_{n}(R,\delta,W)$\ we obtain that $\varphi(\delta
(s))W=\delta_{n}(\varphi(s))W$

\noindent$=W\varphi(s)$. We proceed by induction and assume that
$\varphi(\delta^{i}(s))W^{i}=W^{i}\varphi(s)$\ holds for some $i\geq1$.\ Now
the substitution of $\delta(s)$ into the place of $s$\ gives%
\[
\varphi(\delta^{i+1}(s))W^{i+1}=\varphi(\delta^{i}(\delta(s)))W^{i}%
W=W^{i}\varphi(\delta(s))W=W^{i}W\varphi(s)=W^{i+1}\varphi(s).
\]
Thus we have%
\[
\varphi^{(\delta)}((rw^{i})(sw^{j}))=\varphi^{(\delta)}(r\delta^{i}%
(s)w^{i+j})=\varphi(r\delta^{i}(s))W^{i+j}z^{i+j}=
\]%
\[
\varphi(r)\varphi(\delta^{i}(s))W^{i}W^{j}z^{i+j}=\varphi(r)W^{i}%
\varphi(s)W^{j}z^{i+j}=
\]%
\[
\varphi(r)W^{i}z^{i}\varphi(s)W^{j}z^{j}=\varphi^{(\delta)}(rw^{i}%
)\varphi^{(\delta)}(sw^{j}).
\]

\noindent(1) It follows from $\varphi(cr_{i})=c\varphi(r_{i})$.

\noindent(2) It follows from the fact, that $\varphi(r_{i})W^{i}=0$ implies
$r_{i}=0$.

\noindent(3) The application of $\delta_{n}\circ\varphi=\varphi\circ\delta$
and $\delta_{n}(W)=W$ gives%
\[
(\delta_{z})_{n}(\varphi^{(\delta)}(w))W=(\delta_{z})_{n}(\varphi
(r)Wz)W=\delta_{n}(\varphi(r)W)zW=
\]%
\[
\delta_{n}(\varphi(r))\delta_{n}(W)zW=\varphi(\delta(r))WWz=W\varphi
(r)Wz=W\varphi^{(\delta)}(w),
\]
whence $\varphi^{(\delta)}(w)\in\mathrm{M}_{n}(R[z],\delta_{z},W)$ follows.
Since $\varphi^{(\delta)}(r_{i})=\varphi(r_{i})\in\mathrm{M}_{n}%
(R[z],\delta_{z},W)$, we obtain that $\varphi^{(\delta)}(f(w))\in
\mathrm{M}_{n}(R[z],\delta_{z},W)$ for all $f(w)\in R[w,\delta]$. $\square$

\bigskip

\noindent\textbf{2.4. Corollary.}\textit{ For an automorphism }$\delta
\in\mathrm{Aut}_{K}(R)$\textit{ with }$\delta^{n}=\mathrm{id}_{R}$\textit{,
the number of elements of the cyclic subgroup }$\left\langle \delta
\right\rangle =\{\delta^{i}\mid1\leq i\leq n\}$\textit{ of }$\mathrm{Aut}%
_{K}(R)$\textit{ is a divisor of }$n$\textit{ (and equality not necessarily
holds}). \textit{A natural diagonal embedding }$\Delta:R\longrightarrow
\mathrm{M}_{n}(R)$\textit{ of }$K$\textit{-algebras can be defined by }%
$\Delta(r)=\Sigma_{i=1}^{n}\delta^{i}(r)E_{i,i}$\textit{, where }$r\in
R$\textit{. Using }$H=E_{1,2}+E_{2,3}+\cdots+E_{n-1,n}+E_{n,1}$\textit{, a
straightforward calculation shows that }$\delta_{n}(\Delta(r))H=H\Delta
(r)$\textit{. It follows that }$\Delta$\textit{ is an }$R\longrightarrow
\mathrm{M}_{n}(R,\delta,H)$\textit{ map. Since }$\Delta\circ\delta=\delta
_{n}\circ\Delta$\textit{ and }$H\in\mathrm{GL}_{n}(K)\cap\mathrm{M}%
_{n}(\mathrm{Fix}(\delta))$\textit{, there is a unique homomorphic extension
}$\Delta^{(\delta)}:R[w,\delta]\longrightarrow\mathrm{M}_{n}(R[z])$\textit{ of
}$\Delta$\textit{ such that }$\Delta^{(\delta)}(r)=\Delta(r)$\textit{,
}$\Delta^{(\delta)}(w)=Hz$\textit{ and}%
\[
\Delta^{(\delta)}(r_{0}+r_{1}w+\cdots+r_{t}w^{t})=\Delta(r_{0})+\Delta
(r_{1})Hz+\cdots+\Delta(r_{t})H^{t}z^{t}.
\]
\textit{This map is a}%
\[
\Delta^{(\delta)}:R[w,\delta]\longrightarrow\mathrm{M}_{n}(R[z],\delta_{z},H)
\]
\textit{embedding of }$K$\textit{-algebras. If }$c\in\mathrm{Fix}(\delta
)$\textit{, then }$\Delta^{(\delta)}(cf(w))=c\Delta^{(\delta)}(f(w))$\textit{
for all }$f(w)\in R[w,\delta]$\textit{. Now }$cH=Hc$\textit{ ensures that
}$\Delta^{(\delta)}(f(w)c)=\Delta^{(\delta)}(f(w))c$\textit{.}

\noindent\textit{If }$\left\vert \left\langle \delta\right\rangle \right\vert
=n$\textit{, then the choice }$G=\{\sigma_{1}=\delta,\sigma_{2}=\delta
^{2},\ldots,\sigma_{n}=\delta^{n}\}$\textit{ in Theorem 2.2 gives that
}$\Delta=\Gamma_{G}$\textit{ and }$H=P_{\delta}$\textit{.}

\bigskip

\noindent3. THE\ LIE NILPOTENT\ CAYLEY-HAMILTON THEOREM

\bigskip

A Lie nilpotent analogue of classical determinant theory was developed in
[S1], further details can be found in [Do,SvW]. Here we present the basic
definitions and results about the sequences of right determinants and right
characteristic polynomials, including the so-called Lie nilpotent right
Cayley-Hamilton identities.

For an $n\times n$ matrix $A=[a_{i,j}]$ over an arbitrary (possibly
non-commutative) ring or algebra $R$ with $1$, the element%
\[
\mathrm{sdet}(A)=\underset{\tau,\pi\in\mathrm{S}_{n}}{\sum}\mathrm{sgn}%
(\pi)a_{\tau(1),\pi(\tau(1))}\cdots a_{\tau(t),\pi(\tau(t))}\cdots
a_{\tau(n),\pi(\tau(n))}%
\]
of $R$ is called the symmetric determinant of $A$. The symmetric adjoint
$A^{\ast}=[a_{r,s}^{\ast}]$ of $A=[a_{i,j}]$ is defined as the following
natural symmetrization of the classical adjoint:%
\[
a_{r,s}^{\ast}=\underset{\tau,\pi}{\sum}\mathrm{sgn}(\pi)a_{\tau(1),\pi
(\tau(1))}\cdots a_{\tau(s-1),\pi(\tau(s-1))}a_{\tau(s+1),\pi(\tau
(s+1))}\cdots a_{\tau(n),\pi(\tau(n))}%
\]
where the sum is taken over all $\tau,\pi\in\mathrm{S}_{n}$ with $\tau(s)=s$
and $\pi(s)=r$. We note that the $(r,s)$ entry of $A^{\ast}$ is exactly the
signed symmetric determinant $(-1)^{r+s}\mathrm{sdet}(A_{s,r})$\ of the
$(n-1)\times(n-1)$\ minor $A_{s,r}$\ of $A$ arising from the deletion of the
$s$-th row and the $r$-th column of $A$. If $R$\ is commutative, then
$\mathrm{sdet}(A)=n!\mathrm{\det}(A)$ and $A^{\ast}=(n-1)!\mathrm{adj}(A)$,
where $\mathrm{\det}(A)$ and $\mathrm{adj}(A)$ denote the ordinary determinant
and adjoint of $A$.

\noindent The next result of Domokos plays a fundamental role in the proof of
Theorem 4.1, on which the rest of section (4) is based.

\bigskip

\noindent\textbf{3.1. Theorem ([Do]]). }\textit{Let }$R$\textit{ be an algebra
over a field }$K$\textit{ of characteristic zero. If }$A\in\mathrm{M}_{n}%
(R)$\textit{ and }$T\in\mathrm{GL}_{n}(K)$\textit{, then }$(T^{-1}AT)^{\ast
}=T^{-1}A^{\ast}T$\textit{.}

\bigskip

The right adjoint sequence $(P_{k})_{k\geq1}$ of a matrix $A\in\mathrm{M}%
_{n}(R)$ is defined by the recursion: $P_{1}=A^{\ast}$ and $P_{k+1}%
=(AP_{1}\cdots P_{k})^{\ast}$ for $k\geq1$. It is easy to see that for any
$i\geq1$, the right adjoint sequence of the matrix $AP_{1}\cdots P_{i}$ is
$(P_{k})_{k\geq i+1}$. The $k$-th right adjoint of $A$ is defined as%
\[
\mathrm{radj}_{(k)}(A)=nP_{1}\cdots P_{k}.
\]

The $k$-th right determinant of $A$ is the trace of $AP_{1}\cdots P_{k}$:%
\[
\mathrm{rdet}_{(k)}(A)=\mathrm{tr}(AP_{1}\cdots P_{k}).
\]
We note that
\[
\mathrm{rdet}_{(1)}(A)=\mathrm{tr}(AA^{\ast})=\mathrm{sdet}(A)=\mathrm{tr}%
(A^{\ast}A)
\]
and the following theorem shows that $\mathrm{radj}_{(k)}(A)$ and
$\mathrm{rdet}_{(k)}(A)$ can play a role similar to that played by the
ordinary adjoint and determinant in the commutative case.

\bigskip

\noindent\textbf{3.2. Theorem ([S1],[SvW]).}\textit{ If }$\frac{1}{n}\in
R$\textit{ and the ring }$R$\textit{ is Lie nilpotent of index }$k$\textit{,
then for a matrix }$A\in\mathrm{M}_{n}(R)$\textit{\ we have}%
\[
A\mathrm{radj}_{(k)}(A)=nAP_{1}\cdots P_{k}=\mathrm{rdet}_{(k)}(A)I_{n}.
\]

\bigskip

The above Theorem 3.2 is not used explicitly, however it helps our
understanding and serves as a starting point in the proof of Theorem 3.4.

Let $R[x]$ denote the ring of polynomials of the single commuting
indeterminate $x$, with coefficients in $R$. The $k$-th right characteristic
polynomial of $A$ is the $k$-th right determinant of the $n\times n$ matrix
$xI_{n}-A$ in $\mathrm{M}_{n}(R[x])$:%
\[
p_{A,k}(x)=\mathrm{rdet}_{(k)}(xI_{n}-A).
\]
\noindent\textbf{3.3. Proposition ([SvW]).}\textit{ The }$k$\textit{-th right
characteristic polynomial}

\noindent$p_{A,k}(x)\in R[x]$\textit{ of }$A\in\mathrm{M}_{n}(R)$\textit{\ is
of the form}%
\[
p_{A,k}(x)=\lambda_{0}^{(k)}+\lambda_{1}^{(k)}x+\cdots+\lambda_{n^{k}-1}%
^{(k)}x^{n^{k}-1}+\lambda_{n^{k}}^{(k)}x^{n^{k}},
\]
\textit{where }$\lambda_{0}^{(k)},\lambda_{1}^{(k)},\ldots,\lambda_{n^{k}%
-1}^{(k)},\lambda_{n^{k}}^{(k)}\in R$\textit{ and }$\lambda_{n^{k}}%
^{(k)}=n\left\{  (n-1)!\right\}  ^{1+n+n^{2}+\cdots+n^{k-1}}$\textit{.}

\bigskip

The degree and the leading coefficient of the $k$-th right characteristic
polynomial in Proposition 3.3 will provide the degree of right integrality in
Theorems (A)\ and~(B).

\bigskip

\noindent\textbf{3.4. Theorem ([S1],[SvW]).}\textit{ If }$\frac{1}{n}\in
R$\textit{ and the ring }$R$\textit{ is Lie nilpotent of index }$k$\textit{,
then a right Cayley-Hamilton identity}%
\[
(A)p_{A,k}=I_{n}\lambda_{0}^{(k)}+A\lambda_{1}^{(k)}+\cdots+A^{n^{k}-1}%
\lambda_{n^{k}-1}^{(k)}+A^{n^{k}}\lambda_{n^{k}}^{(k)}=0
\]
\textit{with right scalar coefficients holds for }$A\in\mathrm{M}_{n}%
(R)$\textit{. We also have }$(A)u=0$\textit{, where }$u(x)=p_{A,k}%
(x)h(x)$\textit{ and }$h(x)\in R[x]$\textit{ is arbitrary.}

\bigskip

\noindent4. THE\ RIGHT CHARACTERISTIC\ POLYNOMIALS\ OF

\noindent\ \ \ \ A\ SKEW CENTRALIZING\ MATRIX

\bigskip

\noindent\textbf{4.1. Theorem.}\textit{ Let }$R$\textit{ be an algebra over a
field }$K$\textit{ of characteristic zero,}

\noindent$\delta:R\longrightarrow R$\textit{ be a }$K$\textit{-endomorphism
and }$W\in\mathrm{GL}_{n}(K)$\textit{ be an invertible matrix. If }%
$A\in\mathrm{M}_{n}(R,\delta,W)$\textit{ is a }$(\delta,W)$%
\textit{-centralizing matrix, then }$A^{\ast}\in\mathrm{M}_{n}(R,\delta
,W)$\textit{. In other words, the matrix algebra }$\mathrm{M}_{n}(R,\delta
,W)$\textit{\ is closed with respect to taking the symmetric adjoint.}

\bigskip

\noindent\textbf{Proof.} The definition of $A^{\ast}$ and the element-wise
action of $\delta_{n}$ ensure that $\delta_{n}(A^{\ast})=\delta_{n}(A)^{\ast}$
for all $A\in\mathrm{M}_{n}(R)$. Since $W\in\mathrm{GL}_{n}(K)$, Theorem 3.1
gives that

\noindent$(WAW^{-1})^{\ast}=WA^{\ast}W^{-1}$ for all $A\in\mathrm{M}_{n}(R)$.
Thus for a matrix $A\in\mathrm{M}_{n}(\!R,\delta,W)$\ we have%
\[
\delta_{n}(A^{\ast})=\delta_{n}(A)^{\ast}=(WAW^{-1})^{\ast}=WA^{\ast}W^{-1},
\]
whence $A^{\ast}\in\mathrm{M}_{n}(R,\delta,W)$ follows. $\square$

\bigskip

\noindent\textbf{Remark.} If $A\in\mathrm{M}_{n,d}(E)=\mathrm{M}%
_{n}(E,\varepsilon,Q_{d})$ as in Example (1), then the application of Theorem
4.1 gives that $A^{\ast}\in\mathrm{M}_{n,d}(E)=\mathrm{M}_{n}(E,\varepsilon
,Q_{d})$. This is one of the main results in [S2].

\newpage

\noindent\textbf{4.2. Theorem.}\textit{ Let }$R$\textit{ be an algebra over a
field }$K$\textit{ of characteristic zero,}

\noindent$\delta:R\longrightarrow R$\textit{ be a }$K$\textit{-endomorphism
and }$W\in\mathrm{GL}_{n}(K)$\textit{ be an invertible matrix. If }%
$A\in\mathrm{M}_{n}(R,\delta,W)$\textit{ is a }$(\delta,W)$%
\textit{-centralizing matrix and }$k\geq1$\textit{ is an integer, then for the
}$k$\textit{-th right determinant we have }$\mathrm{rdet}_{(k)}(A)\in
\mathrm{Fix}(\delta)$\textit{.}

\bigskip

\noindent\textbf{Proof.} The repeated application of Theorem 4.1 gives that
the recursion

\noindent$P_{1}=A^{\ast}$ and $P_{k+1}=(AP_{1}\cdots P_{k})^{\ast}$ starting
from a $(\delta,W)$-centralizing matrix $A\in\mathrm{M}_{n}(R,\delta,W)$
defines a sequence $(P_{k})_{k\geq1}$ in $\mathrm{M}_{n}(R,\delta,W)$. Since
$\mathrm{rdet}_{(k)}(A)=\mathrm{tr}(AP_{1}\cdots P_{k})$ is the trace of
$AP_{1}\cdots P_{k}\in\mathrm{M}_{n}(R,\delta,W)$ and the trace of a
$(\delta,W)$\textit{-}centralizing matrix is in $\mathrm{Fix}(\delta)$ by
Proposition 2.1, the proof is complete. $\square$

\bigskip

\noindent\textbf{4.3. Corollary.}\textit{ Let }$R$\textit{ be an algebra over
a field }$K$\textit{ of characteristic zero,}

\noindent$\delta:R\longrightarrow R$\textit{ be a }$K$\textit{-endomorphism
and }$W\in\mathrm{GL}_{n}(K)$\textit{ be an invertible matrix}. \textit{If
}$A\in\mathrm{M}_{n}(R,\delta,W)$\textit{ is a }$(\delta,W)$%
\textit{-centralizing matrix and }$k\geq1$\textit{ is an integer, then for the
}$k$\textit{-th right characteristic polynomial we have }$p_{A,k}%
(x)\in\mathrm{Fix}(\delta)[x]$\textit{.} \textit{In other words, the
coefficients of the }$k$\textit{-th right }$p_{A,k}(x)=\mathrm{rdet}%
_{(k)}(xI_{n}-A)$\textit{ characteristic polynomial are in }$\mathrm{Fix}%
(\delta)$\textit{.}

\bigskip

\noindent\textbf{Proof.} We use the natural extension $\delta_{x}%
:R[x]\longrightarrow R[x]$\ of $\delta$. Since%
\[
xI_{n}-A\in\mathrm{M}_{n}(R[x],\delta_{x},W),
\]
Theorem 4.2 gives that $p_{A,k}(x)=\mathrm{rdet}_{(k)}(xI_{n}-A)$ is in
$\mathrm{Fix}(\delta_{x})=\mathrm{Fix}(\delta)[x]$. $\square$

\bigskip

\noindent\textbf{Remark.} If $A\in\mathrm{M}_{n,d}(E)=\mathrm{M}%
_{n}(E,\varepsilon,Q_{d})$ as in Example (1), then the application of
Corollary 4.3 gives that $p_{A,2}(x)\in\mathrm{Fix}(\varepsilon)[x]=E_{0}[x]$.
This fact was exploited in [S2].

\bigskip

\noindent\textbf{4.4. Corollary.}\textit{ Let }$R=R_{0}\oplus R_{1}%
\oplus\cdots\oplus R_{n-1}$\textit{ be a }$\mathbb{Z}_{n}$\textit{-graded
algebra over a field }$K$\textit{ of characteristic zero and }$e\in K$\textit{
be a primitive }$n$\textit{-th root of unity. If }$k\geq1$\textit{ is an
integer and }$A\in\mathrm{M}_{n}^{g}(R)$\textit{ is a "graded" }$n\times
n$\textit{ matrix with respect to the \ given }$\mathbb{Z}_{n}$%
\textit{-grading, then for the }$k$\textit{-th right determinant and }%
$k$\textit{-th right characteristic polynomial of }$A$\textit{ we have}%
\[
\mathrm{rdet}_{(k)}(A)\in R_{0}\text{\textit{ and }}p_{A,k}(x)=\mathrm{rdet}%
_{(k)}(xI_{n}-A)\in R_{0}[x].
\]
\textbf{Proof.} Since $\mathrm{M}_{n}^{g}(R)=\mathrm{M}_{n}(R,\widehat
{e},D_{e})$ (see Example (2)) and $\mathrm{Fix}(\widehat{e})=R_{0}$, Theorem
4.2 and Corollary 4.3 can be applied. $\square$

\bigskip

\noindent\textbf{Remark.} In the absence of a primitive $n$-th root of unity,
the direct proof (not using Example (2), Theorem 4.2 and Corollary 4.3) of the
implications%
\[
A\in\mathrm{M}_{n}^{g}(R)\Longrightarrow A^{\ast}\in\mathrm{M}_{n}%
^{g}(R),\text{ }\mathrm{rdet}_{(k)}(A)\in R_{0},\text{ }p_{A,k}%
(x)=\mathrm{rdet}_{(k)}(xI_{n}-A)\in R_{0}[x]
\]
is rather technical.

\bigskip

\noindent The combination of Theorem 3.4 and Corollary 4.3 gives the following.

\bigskip

\noindent\textbf{4.5. Theorem.}\textit{ Let }$R$\textit{ be an algebra over a
field }$K$\textit{ of characteristic zero,}

\noindent$\delta:R\longrightarrow R$\textit{ be an endomorphism and }%
$W\in\mathrm{GL}_{n}(K)$\textit{ be an invertible matrix}. \textit{If }%
$R$\textit{ is Lie nilpotent of index }$k$\textit{\ and }$A\in\mathrm{M}%
_{n}(R,\delta,W)$\textit{, then a right Cayley-Hamilton identity}%
\[
(A)p_{A,k}=I_{n}\lambda_{0}^{(k)}+A\lambda_{1}^{(k)}+\cdots+A^{n^{k}-1}%
\lambda_{n^{k}-1}^{(k)}+A^{n^{k}}\lambda_{n^{k}}^{(k)}=0
\]
\textit{holds, where the coefficients }$\lambda_{i}^{(k)}$\textit{, }$0\leq
i\leq n^{k}$\textit{\ of }$p_{A,k}(x)=\mathrm{rdet}_{(k)}(xI_{n}-A)$\textit{
are in }$\mathrm{Fix}(\delta)$\textit{. Since }$\lambda_{n^{k}}^{(k)}%
=n\left\{  (n-1)!\right\}  ^{1+n+n^{2}+\cdots+n^{k-1}}$\textit{ is invertible
(in }$K$\textit{), the above identity provides the right integrality of
}$\mathrm{M}_{n}(R,\delta,W)$\textit{ over }$\mathrm{Fix}(\delta)$\textit{ of
degree }$n^{k}$\textit{.}

\bigskip

\noindent5. THE\ PROOFS\ OF\ THEOREM\ (A)\ AND\ (B)

\bigskip

\noindent\textbf{Proof of Theorem (A)}. Take an element $\tau\in G$\ and
consider the embedding $\Gamma_{G}:\!R\!\longrightarrow\!\mathrm{M}_{n}%
(R,\tau,P_{\tau})$ in Theorem 2.2. For $r\in R$ Theorem 4.5 ensures that
\[
I_{n}\lambda_{0}^{(k)}+\Gamma_{G}(r)\lambda_{1}^{(k)}+\cdots+(\Gamma
_{G}(r))^{n^{k}-1}\lambda_{n^{k}-1}^{(k)}+(\Gamma_{G}(r))^{n^{k}}%
\lambda_{n^{k}}^{(k)}=0
\]
holds for the $(\tau,P_{\tau})$-centralizing matrix $\Gamma_{G}(r)\in
\mathrm{M}_{n}(R,\tau,P_{\tau})$, where%
\[
p_{\Gamma_{G}(r),k}(x)=\mathrm{rdet}_{(k)}(xI_{n}-\Gamma_{G}(r))=\lambda
_{0}^{(k)}+\lambda_{1}^{(k)}x+\cdots+\lambda_{n^{k}-1}^{(k)}x^{n^{k}%
-1}+\lambda_{n^{k}}^{(k)}x^{n^{k}},
\]
$\lambda_{n^{k}}^{(k)}=q\geq1$ is an integer and $\lambda_{i}^{(k)}%
\in\mathrm{Fix}(\tau)$ for each $0\leq i\leq n^{k}$. Thus we have

\noindent$c_{i}=\frac{1}{q}\lambda_{i}^{(k)}\in\mathrm{Fix}(G)$ for each index
$0\leq i\leq n^{k}-1$. It follows that%
\[
\Gamma_{G}(c_{0}+rc_{1}+\cdots+r^{n^{k}-1}c_{n^{k}-1}+r^{n^{k}})=
\]%
\[
I_{n}c_{0}+\Gamma_{G}(r)c_{1}+\cdots+(\Gamma_{G}(r))^{n^{k}-1}c_{n^{k}%
-1}+(\Gamma_{G}(r))^{n^{k}}=
\]%
\[
\left(  I_{n}\lambda_{0}^{(k)}+\Gamma_{G}(r)\lambda_{1}^{(k)}+\cdots
+(\Gamma_{G}(r))^{n^{k}-1}\lambda_{n^{k}-1}^{(k)}+(\Gamma_{G}(r))^{n^{k}%
}\lambda_{n^{k}}^{(k)}\right)  \cdot\frac{1}{q}=0
\]
and $\ker(\Gamma_{G})=\{0\}$ gives the desired right integrality. $\square$

\bigskip

\noindent\textbf{Proof of Theorem (B)}. Consider the embeddings%
\[
\Delta:R\longrightarrow\mathrm{M}_{n}(R,\delta,H)\text{ and }\Delta^{(\delta
)}:R[w,\delta]\longrightarrow\mathrm{M}_{n}(R[z],\delta_{z},H)
\]
in Corollary 2.4, where $\Delta(r)=\Sigma_{i=1}^{n}\delta^{i}(r)E_{i,i}$,
$H=E_{1,2}+E_{2,3}+\cdots+E_{n-1,n}+E_{n,1}$ and%
\[
\Delta^{(\delta)}(r_{0}+r_{1}w+\cdots+r_{t}w^{t})=\Delta(r_{0})+\Delta
(r_{1})Hz+\cdots+\Delta(r_{t})H^{t}z^{t}.
\]
Since $R[z]$ is also Lie nilpotent of index $k$, Theorem 4.5 gives that%
\[
I_{n}\lambda_{0}^{(k)}(z)+\Delta^{(\delta)}(f(w))\lambda_{1}^{(k)}(z)+\cdots
\]%
\[
\cdots+(\Delta^{(\delta)}(f(w))^{n^{k}-1}\lambda_{n^{k}-1}^{(k)}%
(z)+(\Delta^{(\delta)}(f(w))^{n^{k}}\lambda_{n^{k}}^{(k)}(z)=0,
\]
where the coefficients of the $k$-th right characteristic polynomial%
\[
\mathrm{rdet}_{(k)}(xI_{n}-\Delta^{(\delta)}(f(w)))=\lambda_{0}^{(k)}%
(z)+\lambda_{1}^{(k)}(z)x+\cdots+\lambda_{n^{k}-1}^{(k)}(z)x^{n^{k}-1}%
+\lambda_{n^{k}}^{(k)}(z)x^{n^{k}}%
\]
of the $(\delta_{z},H)$-centralizing matrix $\Delta^{(\delta)}(f(w))\in
\mathrm{M}_{n}(R[z],\delta_{z},H)$ are in $\mathrm{Fix}(\delta_{z}%
)=\mathrm{Fix}(\delta)[z]$ and $\lambda_{n^{k}}^{(k)}(z)=q\geq1$ is an integer.

\noindent According to Example (2), the natural $\mathbb{Z}_{n}$-grading%
\[
R[z]=R[z^{n}]\oplus zR[z^{n}]\oplus\cdots\oplus z^{n-1}R[z^{n}]
\]
of the polynomial algebra $R[z]$\ defines $\mathrm{M}_{n}^{g}(R[z])$ as%
\[
\mathrm{M}_{n}^{g}(R[z])=\{A=[a_{i,j}(z)]\in\mathrm{M}_{n}(R[z])\mid
a_{i,j}(z)\in z^{j-i}R[z^{n}]\text{ for all }1\leq i,j\leq n\}=
\]%
\[
\left[
\begin{array}
[c]{ccccc}%
R[z^{n}] & zR[z^{n}] & \cdots & z^{n-2}R[z^{n}] & z^{n-1}R[z^{n}]\\
z^{n-1}R[z^{n}] & R[z^{n}] & zR[z^{n}] & \ddots & z^{n-2}R[z^{n}]\\
\vdots & z^{n-1}R[z^{n}] & \ddots & \ddots & \vdots\\
z^{2}R[z^{n}] & \ddots & \ddots & R[z^{n}] & zR[z^{n}]\\
zR[z^{n}] & z^{2}R[z^{n}] & \cdots & z^{n-1}R[z^{n}] & R[z^{n}]
\end{array}
\right]  =\mathrm{M}_{n}(R[z],\widehat{e},D_{e}).
\]
Clearly, $\Delta(r_{i})\!\in\!\mathrm{M}_{n}^{g}(R[z])$ and $Hz\!\in
\!\mathrm{M}_{n}^{g}(R[z])$ imply that $\Delta(r_{i})H^{i}z^{i}\!\in
\!\mathrm{M}_{n}^{g}(R[z])$ for all $i\geq0$. Thus we have $\Delta^{(\delta
)}(f(w))\in\mathrm{M}_{n}^{g}(R[z])$ for all $f(w)\in R[w,\delta]$. The use of
Corollary 4.4 gives that
\[
\mathrm{rdet}_{(k)}(xI_{n}-\Delta^{(\delta)}(f(w)))\in(R[z^{n}])[x],
\]
whence $\lambda_{i}^{(k)}(z)=\mu_{i}(z^{n})\in\mathrm{Fix}(\delta)[z^{n}]$
follows for all $0\leq i\leq n^{k}-1$.

\noindent Since $H^{n}=I_{n}$, for any polynomial $\mu(z^{n})\in
\mathrm{Fix}(\delta)[z^{n}]$ there exists a (unique) polynomial $g(w^{n}%
)\in\mathrm{Fix}(\delta)[w^{n}]$ such that%
\[
\mu(z^{n})I_{n}=c_{0}I_{n}+c_{1}z^{n}I_{n}+\cdots+c_{d}z^{nd}I_{n}=
\]%
\[
\Delta(c_{0})+\Delta(c_{1})H^{n}z^{n}+\cdots+\Delta(c_{d})H^{nd}z^{nd}=
\]%
\[
\Delta^{(\delta)}(c_{0}+c_{1}w^{n}+\cdots+c_{d}w^{nd})=\Delta^{(\delta
)}(g(w^{n})),
\]
where the coefficients $c_{0},c_{1},\ldots,c_{d}\in\mathrm{Fix}(\delta)$ of
$\mu(z^{n})$ and $g(w^{n})$ coincide. For each $0\leq i\leq n^{k}-1$ take
$g_{i}(w^{n})\in\mathrm{Fix}(\delta)[w^{n}]$ such that%
\[
\Delta^{(\delta)}(g_{i}(w^{n}))=\frac{1}{q}\mu_{i}(z^{n})I_{n}=\frac{1}%
{q}\lambda_{i}^{(k)}(z)I_{n}.
\]
Thus%
\[
\Delta^{(\delta)}\left(  g_{0}(w^{n})+f(w)g_{1}(w^{n})+\cdots+f^{n^{k}%
-1}(w)g_{n^{k}-1}(w^{n})+f^{n^{k}}(w)\right)  =
\]%
\[
\frac{1}{q}\left(  I_{n}\lambda_{0}^{(k)}(z)+\Delta^{(\delta)}(f(w))\lambda
_{1}^{(k)}(z)+\cdots\right.
\]%
\[
\left.  \cdots+(\Delta^{(\delta)}(f(w))^{n^{k}-1}\lambda_{n^{k}-1}%
^{(k)}(z)+(\Delta^{(\delta)}(f(w))^{n^{k}}\lambda_{n^{k}}^{(k)}(z)\right)  =0
\]
and $\ker(\Delta^{(\delta)})=\{0\}$ gives the desired integrality. $\square$

\bigskip

\noindent REFERENCES

\bigskip

\noindent\lbrack Do] M. Domokos, \textit{Cayley-Hamilton theorem for }%
$2\times2$\textit{\ matrices over the Grassmann algebra}, J. Pure Appl.
Algebra 133 (1998), 69-81.

\noindent\lbrack Ke] A. R. Kemer,\textit{\ Ideals of Identities of Associative
Algebras,} Translations of Math. Monographs, Vol. 87 (1991), AMS, Providence,
Rhode Island.

\noindent\lbrack S1] J. Szigeti, \textit{New determinants and the
Cayley-Hamilton theorem for matrices over Lie nilpotent rings}, Proc. Amer.
Math. Soc. 125 (1997), 2245-2254.

\noindent\lbrack S2] J. Szigeti, \textit{On the characteristic polynomial of
supermatrices}, Israel Journal of Mathematics Vol. 107 (1998), 229-235.

\noindent\lbrack SvW] J. Szigeti and L. van Wyk,\textit{ Determinants for
}$n\times n$\textit{ matrices and the symmetric Newton formula in the
}$3\times3$\textit{ case,} Linear and Multilinear Algebra, Vol. 62, No. 8
(2014), 1076-1090.

\end{document}